\documentclass[]{amsart}
\usepackage[foot]{amsaddr}
\usepackage[a4paper,top=2.5cm,bottom=3cm,inner=3cm,outer=3cm]{geometry}

\usepackage[utf8]{inputenc}

\usepackage{amsmath,amsthm}
\usepackage{amsfonts,amssymb,mathrsfs}

\usepackage{xy}
\usepackage[bookmarks=false,breaklinks=true]{hyperref}
\usepackage{float}
\usepackage{graphicx}
\usepackage[percent]{overpic}
\usepackage{color}

\usepackage{epigraph}
\definecolor{myurlcolor}{rgb}{0,0,0.7}
\usepackage{hyperref}
\hypersetup{colorlinks,
linkcolor=myurlcolor,
citecolor=myurlcolor,
urlcolor=myurlcolor}
\usepackage[mathscr]{euscript}

\input xy
\xyoption{all}
\newcommand{\maps}{\colon}    %correct symbol for colon in f: X -> Y
\newcommand{\Z}{{\mathbb Z}}  %integers
\newcommand{\R}{{\mathbb R}}  %real numbers
\newcommand{\define}[1]{{\bf \boldmath{#1}}}

\theoremstyle{definition}

        \newcommand{\be}{\begin{equation}}
        \newcommand{\ee}{\end{equation}}
        \newcommand{\ba}{\begin{eqnarray}}
        \newcommand{\ea}{\end{eqnarray}}
        \newcommand{\ban}{\begin{eqnarray*}}
        \newcommand{\ean}{\end{eqnarray*}}
        \newcommand{\barr}{\begin{array}}
        \newcommand{\earr}{\end{array}}

\begin{document}
\title{The Kuramoto--Sivashinky Equation}
\author[Baez]{John C.\ Baez} 
\address{Department of Mathematics, University of California, Riverside CA, 92521, USA}
\address{Centre for Quantum Technologies, National University of Singapore, 117543, Singapore}
\email{baez@math.ucr.edu}
\author[Huntsman]{Steve Huntsman}
\address{2500 Valley Drive, Alexandria, Virginia, 22302, USA}
\email{sch213@nyu.edu}
\author[Weis]{Cheyne Weis}
\address{James Franck Institute and Department of Physics, University of Chicago, Chicago IL, 60637, USA}
\email{cheyne42@uchicago.edu}
\date{April 30, 2021}
\maketitle

The Kuramoto--Sivashinsky equation 
\[    u_t = -u_{xx} - u_{xxxx} - u_x u \]
applies to a real-valued function $u$ of time $t \in \R$ and space $x \in \R$.  This equation was introduced as a simple 1-dimensional model of instabilities in flames, but it turned out to mathematically fascinating in its own right \cite{Encyclopedia}.  One reason is that the Kuramoto--Sivashinsky equation is a simple model of \emph{Galilean-invariant chaos with an arrow of time}.  

We say this equation is `Galilean invariant' because the Galiei group, the usual group of symmetries in Newtonian mechanics, acts on the set of its solutions.   When space is 1-dimensional, this group is generated by translations in $t$ and $x$, reflections in $x$, and \define{Galilei boosts}, which are transformations to moving coordinate systems:
\[   (t,x) \mapsto (t,x-tv). \]
Translations act in the obvious way.  Spatial reflections act as follows: if $u(t,x)$ is a solution, so is $-u(t,-x)$. Galilei boosts act in a more subtle way: if $u(t,x)$ is a solution, so is $u(t,x-tv) + v$. 

\begin{center}
\includegraphics[width = 45em]{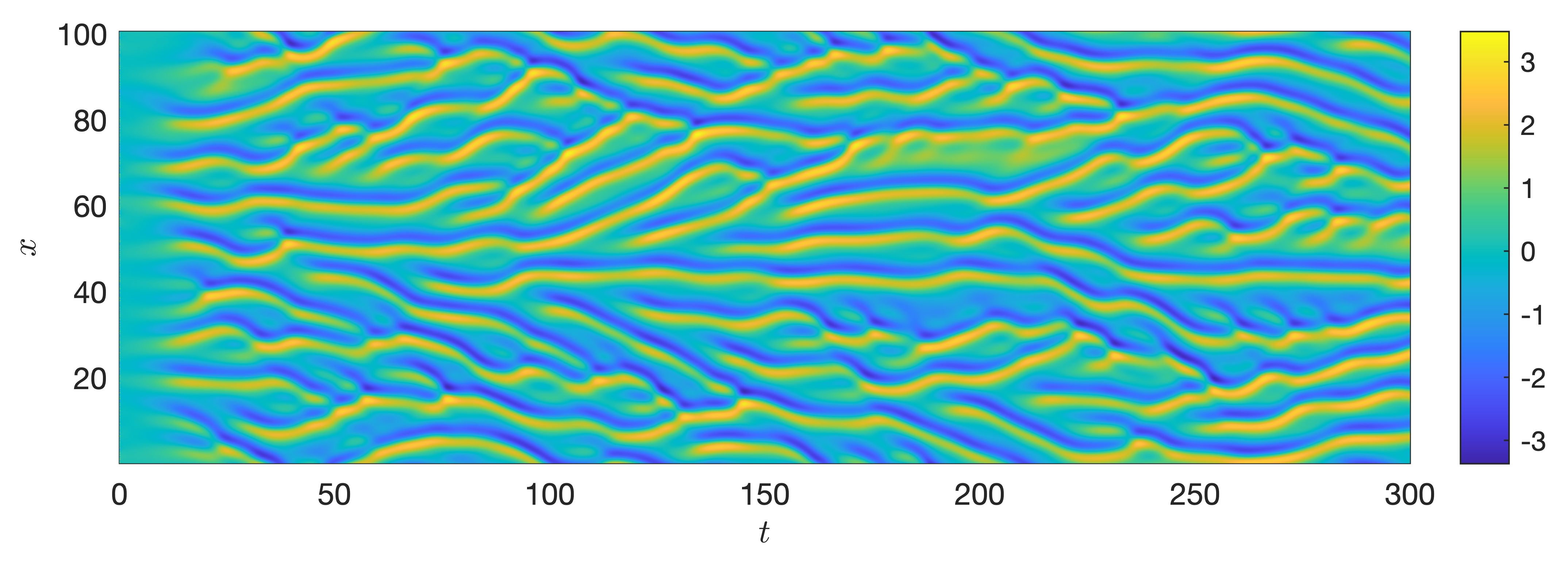} 

Figure 1 --- A solution $u(t,x)$ of the Kuramoto--Sivashinsky equation.  The variable $x$ ranges over the interval $[0,32\pi]$ with its endpoints identified.  Initial data are independent identically distributed random variables, one at each grid point, uniformly distributed in $[-1,1]$.
\end{center}

\vskip 1em

We say the Kuramoto--Sivashinsky equation is `chaotic' because the distance between nearby solutions, defined in a suitable way, can grow exponentially, making the long-term behavior of a solution hard to predict in detail \cite{Chaos}.   And finally, we say this equation has an `arrow of time' because time reversal
\[   (t,x) \mapsto (-t,x) \]
is not a symmetry of this equation.  Indeed, in Figure 1 we see that starting from random initial conditions, manifestly time-asymmetric patterns emerge.  As we move forward in time, it looks as if stripes are born and merge, but never die or split.  Attempting to make this precise leads to an interesting conjecture, but first we need some background.

It is common to study solutions of the Kuramoto--Sivashinsky equations that are spatially periodic, so that $u(t,x) = u(t,x+L)$ for some $L$.  We can then treat space as a circle, the interval $[0,L]$ with its endpoints identified.   For these spatially periodic solutions, the integral $\int_0^L u(t,x) \, dx$ does not change with time.   Applying a Galilean transform adds a constant to this integral.  In what follows we restrict attention to solutions where this integral is zero.  These are roughly the solutions where the stripes are at rest, on average.

We can learn a surprising amount about these solutions by looking at the linearized equation
\[  u_t = -u_{xx} - u_{xxxx} .\]
We can solve this using a Fourier series
\[   u(t,x) = \sum_{0 \ne n \in \Z} \hat{u}_n(t) \, e^{ik_n x} \]
where the frequency of the $n$th mode is $k_n = 2\pi n/L$.   We obtain
\[    \hat{u}_n(t) = \exp\left((k_n^2 - k_n^4) t\right) \, \hat{u}_n(0) .\]
%\[   \frac{d}{dt} \hat{u}_n(t) = (k_n^2 - k_n^4) \hat{u}_n(t) + \sum_{0 \ne m \in \Z} k_m \hat{u}_m(t) \hat{u}_{n - m}(t) .\]
Thus the $n$th mode grows exponentially with time if and only if $k_n^2 - k_n^4 > 0$, which happens when $0 < |n| < 2 \pi L$.   As we increase $L$, more and more modes grow exponentially.  These appear to be the cause of chaos even in the nonlinear equation.  Indeed, all solutions of the Kuramoto--Sivashinsky equation approach an attractive fixed point if $L$ is small enough, but as we increase $L$ we see increasingly complicated behavior, and a `transition to chaos via period doubling', which has been analyzed in great detail \cite{Chaos}.  Interestingly, the nonlinear term stabilizes the exponentially growing modes: in the language of physics, it tends to transfer power from these modes to high-frequency modes, which decay exponentially.

Proving this last fact is not easy.   However, in 1992, Collet, Eckmann, Epstein and Stubbe \cite{ColletAttractor} did this in the process of showing that for any initial data in the Hilbert space
\[   \dot{L}^2 = \{ u \maps 
[0,L] \to \R: \; \int_0^L |u(x)|^2 \, dx < \infty, \; \int_0^L u(x) \, dx = 0 \} \]
the Kuramoto--Sivashinsky equation has a unique solution for $t \ge 0$, in a suitable sense, and that the norm of this solution eventually becomes less than some constant times $L^{8/5}$. These authors also showed any such solution eventually becomes infinitely differentiable, even analytic \cite{ColletAnalyticity}. 

Shortly after this, Temam and Wang went further \cite{TemamWang}.  They showed that all solutions of the Kuramoto--Sivashinsky equation with initial data in $\dot{L}^2$ approach a \emph{finite-dimensional submanifold} of $\dot{L}^2$ as $t \to +\infty$.  They also showed the dimension of this manifold is bounded by a constant times $(\ln L)^{0.2} L^{1.64}$.   This manifold, called the \define{inertial manifold}, describes the `eventual behaviors' of solutions.  

Understanding the eventual behavior of solutions of the Kuramoto--Sivashinsky equation remains a huge challenge.  What are the `stripes' in these solutions?  Can we define them in such a way that stripes are born and merge but never split or disappear, at least after a solution has had time to get close enough to the inertial manifold?  

It is important to note that the visually evident stripes in Figure 1 are not regions where $u$ exceeds some constant, nor regions where it is less than some constant.  Instead, as $x$ increases and $(t,x)$ passes through a stripe, $u(t,x)$ first becomes positive and then negative.  Thus, in the middle of a stripe the derivative $u_x(t,x)$ takes a negative value.   One obvious thing to try, then, is to define a stripe to be a region where $u_x(t,x) < c$ for some suitably chosen negative constant $c$.   Unfortunately the boundaries of these regions, where $u_x(t,x) = c$, tend to be very rough.  Thus, this definition gives many small evanescent `stripes', so the conjecture that eventually stripes are born and merge but never die or split would be false with such a definition. 

One way to solve this problem is to smooth $u_x$ by convolving it with a normalized Gaussian function of $x$.  A bit of experimentation suggests using a normalized Gaussian of standard deviation 2.   Thus, let $v$ equal $u_x$ convolved with this Gaussian, and define a `stripe' to be a region where $v < 0$.  For the solution in Figure 1, these stripes are indicated in Figure 2.  One stripe splits around $t = 10$, but after the solution nears the inertial manifold, stripes never die or split---at least in this example.  We conjecture that this is true generically: that is, for all solutions in some open dense subset of the inertial manifold.    Proving this seems quite challenging, but it might be a step toward rigorously relating the Kuramoto--Sivashinsky equation to a model where stochastically moving particles can appear or merge, but never disappear or split \cite{RostKrug}.  

\begin{center}
\includegraphics[width = 45em]{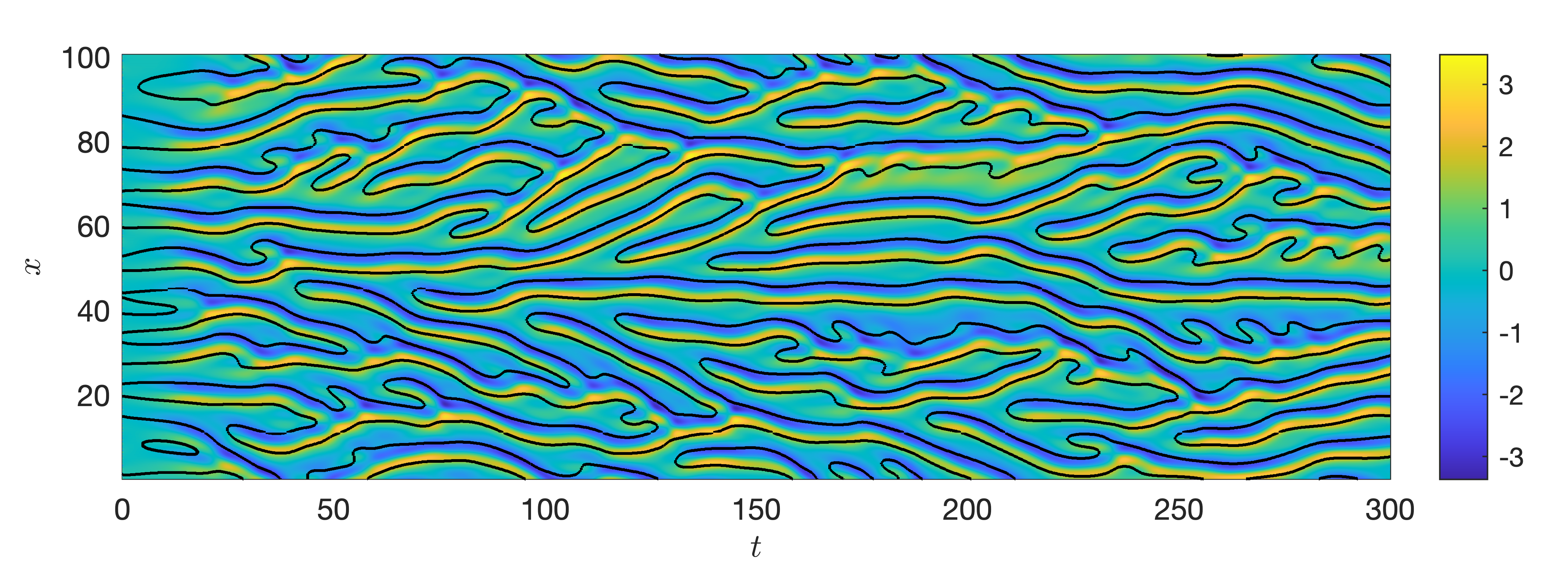} 

Figure 2 --- The same solution of the Kuramoto--Sivashinsky equation, with stripes indicated.
\end{center}
\vskip 1em

Numerical calculations also indicate that generically, solutions eventually have stripes with an average density that approaches about $0.1$ as $L \to +\infty$.  This is close to the inverse of the wavelength of the fastest-growing mode of the linearized equation, $(2^{3/2} \pi)^{-1} \approx 0.1125$.  But we see no clear reason to think these numbers should be exactly equal.    Even proving the existence of a limiting stripe density is an open problem!  Thus, the Kuramoto--Sivashinsky equation continues to pose many mathematical challenges.

\end{document}